\documentclass{article}
\usepackage{amsmath}
\usepackage{amssymb}
 \newtheorem{thm}{Theorem}[section]
 
 \newtheorem{lem}[thm]{Lemma}

 \newtheorem{rem}[thm]{Remark}

 \newtheorem{ex}[thm]{Example}
 \numberwithin{equation}{section}

\newcommand{\de}{{\rm depth}}

\newcommand{\smin}{{\rm Min}}
\newcommand{\ass}{{\rm Ass}}

\newcommand{\cmd}{{\rm cmd}}
    \begin{document}
  \title{Flat local  morphisms of rings with prescribed depth and dimension}
  \author{Cristodor Ionescu\\
  Institute of Mathematics \textit{Simion Stoilow}\\ of the Romanian Academy\\
  P.O. Box 1-764\\ RO 014700 Bucharest  Romania
  \\
  email: cristodor.ionescu@imar.ro}
  \date{}
  \maketitle

  \begin{abstract}
 For any pairs of integers $(n,m)$ and $(d,e)$ such that $0\leq n\leq m,\ 0\leq d\leq e,\ d\leq n,\ e\leq m$ and $n-d\leq m-e$ we construct a local flat ring morphism of noetherian local rings $u:A\to B$ such that $\dim(A)=n, \de(A)=d, \dim(B)=m, \de(B)=e.$
\end{abstract}
\vspace{5mm}

\section{Introduction}
While preparing \cite{Io}, the present author was looking for an example of a flat local ring homomorphism of noetherian local rings $u:(A,m)\to (B,n)$ such that $A$ and $B/mB$ are almost Cohen-Macaulay, while $B$ is not almost Cohen-Macaulay. This means that one should construct for example such a morphism with 
$\de(B)=\de(A)=0, \dim(B)=2$ and $\dim(A)=1.$ Note that actually  the flatness of the homomorphism $u$ is the non-trivial point in the construction. After asking several people without obtaining a satisfactory answer, he decided to let it as an open question  in \cite{Io}. The answer came  soon, an example with the desired features being constructed by Taba\^ a \cite{Ta}. Using his idea, we construct a quite general example of this type, construction that can be useful in various situations. 

\section{The construction}
 We start by pointing out the following easy and well-known fact. 

\begin{lem}\label{exista}
Let k be a field, $n,d\in\mathbb{N}$ such that $0\leq d\leq n.$ Then there exists $m\in\mathbb{N}, m\geq n$ and a monomial ideal $I\subset S:=k[X_1,\ldots,X_m]$ such that 
$\dim(S/I)_{(X_1,\ldots,X_m)}=n$ and $\de(S/I)_{(X_1,\ldots,X_m)}=d.$
\end{lem}

\par\noindent \textbf{\textit{Proof:}} Let $r=n-d.$ If $r=0$ the assertion is clear. Assume that $r>0.$ If $d=0$ let $S=k[X_0,X_1,\ldots, X_r]$  and if $d>0$ let $ S=k[X_0,X_1,\ldots, X_r,T_1,\ldots,T_d].$  Consider for example the monomial ideal $I=(X_0)\cap(X_0,\ldots,X_r)^{r+1}.$ Then 
$\ass(S/I)=\{(X_0),(X_0,\ldots,X_r)\},$ 
hence $\dim(S/I)=r+d=n$ and $\de(S/I)=d.$ Now taking $m=r+1+d$ and renumbering the indeterminates we get the desired relations.

\begin{rem}\label{obser}
Clearly there are also many other choices for a monomial ideal with the properties of the above lemma. For more things about this kind of construction one can see \cite{Sh}.
\end{rem}



\begin{thm}\label{principal}
Let $0\leq d_1\leq n_1$ and $0\leq d_2\leq n_2$ be natural numbers such that $n_1\leq n_2,d_1\leq d_2$ and $n_1-d_1\leq n_2-d_2.$ Then there exists a local flat morphism of noetherian local rings $\ u:(A,m)\to (B,n)$ such that $\de(A)=d_1, \dim(A)=n_1, \de(B)=d_2$ and $\dim(B)=n_2.$
\end{thm}

\par\noindent \textbf{\textit{Proof:}}
Let $k$ be a field and $A=(k[X]/I)_{(X)}, X=(X_1,\ldots,X_m)$ be a local ring obtained cf. \ref{exista} with $\de(A)=d_1, \dim(A)=n_1.$ Let $s=d_2-d_1, t=n_2-n_1$ and let $C=(k[Y]/J)_{(Y)}, Y=(Y_1,\ldots,Y_p)$ be  a local ring obtained cf. \ref{exista} with $\de(C)=s, \dim(C)=t.$ Now let 
$$D:=A\otimes_kC=k[X]/I\otimes_kk[Y]/J=k[X,Y](Ik[X,Y]+Jk[X,Y])$$ and let $B:=D_{(X,Y)}.$ Then obviously the canonical morphism  $u:A\to B$ is flat and local, being a localisation of the base change of the flat morphism $k\to C.$ 
We need the following probably well-known fact:
\begin{lem}\label{minimale}
Let $k$ be a field and  $m,p\in\mathbb{N}.$ Let also $I$ and $J$ be monomial ideals in $k[X]=k[X_1,\ldots,X_m]$ and $k[Y]=k[Y_1,\ldots,Y_p]$ respectively and set $S:=k[X,Y].$ Then $\smin(IS+JS)=\{PS+QS\ \vert P\in\smin(I), Q\in\smin(J)\}.$ Consequently 
$$\dim(S/(IS+JS))=\dim(k[X]/I)+\dim(k[Y]/J).$$
\end{lem}

\par\noindent \textbf{\textit{Proof:}} Using \cite{HM}, 3.4 we obtain 
$$\smin(IS+JS)=\smin(\sqrt{IS+JS})=\smin(\sqrt{\sqrt{IS}+\sqrt{JS}})=$$
$$=\smin(\sqrt{I}S+\sqrt{J}S)=\{PS+QS\ \vert\ P\in\smin(I), Q\in\smin(J)\}.$$

\par\noindent Returning at the proof of the Theorem, by \ref{minimale} we get 
$\dim(B)=\dim(A)+\dim(C)=n_1+t=n_2$ and by \cite{F}, Lemma 2 we have that 
$\de(B)=\de(A)+\de(C)=d_2.$ This concludes the proof of \ref{principal}.

\begin{ex}\label{tabaa}
Let $k$ be a field and let $A=C=(k[X,Y]/(X^2,XY))_{(X,Y)}$ and let
 $B=A\otimes_kC=(k[X,Y,U,V]/(X^2,XY,U^2,UV))_{(X,Y,U,V)}.$ The canonical morphism $u:A\to B$ is the morphism  obtained performing the above construction. This is  the example from \cite{Ta}, namely we have $\dim(A)=1, \dim(B)=2, \de(A)=\de(B)=0.$
\end{ex}

\begin{rem}\label{defect}
If $(A,m)$ is a noetherian local ring, we call the Cohen-Macaulay defect of $A,$ the natural number $\cmd(A)=\dim(A)-\de(A).$ Thus $A$ is Cohen-Macaulay if and only if $\cmd(A)=0$ and $A$ is almost Cohen-Macaulay if and only if $\cmd(A)\leq 1$(see \cite{Io}).
\end{rem}

\begin{ex}\label{fibre} 
Using the above construction, one can also get examples of flat local morphisms of noetherian local rings, whose closed fiber has prescribed Cohen-Macaulay defect, or even more general, has prescribed dimension $n$ and depth $d\leq n.$ Indeed, by (\cite{Mat}, 15.1, 23.3), the flatness of $u$ implies that $n=\dim(B/mB)=n_2-n_1$ and $d=\de(B/mB)=d_2-d_1,$ so that it is enough to choose appropriate values for $n_1\leq  n_2$ and $d_1\leq d_2$ and perform the above construction.
\end{ex}


\vspace{0.4cm}

\par \textbf{Acknowledgment:} The author would like to thank Marius Vl\u adoiu for some illuminating discussions concerning lemma \ref{minimale} and to Javier Majadas for pointing out a missing condition in the statement of the main result.


  \end{document}